\documentclass[letterpaper, 10 pt, conference]{ieeeconf}  

\IEEEoverridecommandlockouts                              
\usepackage[utf8]{inputenc}
\usepackage[T1]{fontenc}
\usepackage{amsmath}
\usepackage{amssymb}
\usepackage{algpseudocode}
\usepackage{amsfonts}
\usepackage{algorithm}
\usepackage{comment}
\usepackage{caption} 
\usepackage{cite}
\usepackage{epsfig}
\usepackage{graphics}   
\usepackage{subcaption}
\usepackage{mathptmx}
\usepackage{multirow}
\usepackage{mathtools}  
\usepackage{times}
\usepackage{xcolor}

\newtheorem{thm}{\textbf{Theorem}}
\newtheorem{lem}{\textbf{Lemma}}

\newtheorem{rem}{\textbf{Remark}}

\newtheorem{definition}{\textbf{Definition}}

\newtheorem{prop}{\textbf{Proposition}}

\title{\LARGE \bf
First and Second Order Optimal $\mathcal{H}_2$ Model Reduction for Linear Continuous-Time Systems
}

\author{Wenshan Zhu, Imad Jaimoukha
\thanks{Wenshan Zhu and I.M.Jaimoukha, are with the Department of Electrical and Electronic Engineering, Imperial College London, London UK.
(wenshan.zhu20@imperial.ac.uk; i.jaimouka@imperial.ac.uk;). }
}

\begin{document}

\maketitle
\thispagestyle{empty}
\pagestyle{empty}

\begin{abstract}
In this paper, we investigate the optimal $\mathcal{H}_2$ model reduction problem for single-input single-output (SISO) continuous-time linear time-invariant (LTI) systems. A semi-definite relaxation (SDR) approach is proposed to determine globally optimal interpolation points, providing an effective way to compute the reduced-order models via Krylov projection-based methods. In contrast to iterative approaches, we use the controllability Gramian and the moment-matching conditions to recast the model reduction problem as a convex optimization by introducing an upper bound $\gamma$ to minimize the $\mathcal{H}_2$ norm of the model reduction error system. We also prove that the relaxation is exact for first order reduced models and demonstrate, through examples, that it is exact for second order reduced models. We compare the performance of our proposed method with other iterative approaches and shift-selection methods on examples. Importantly, our approach also provides a means to verify the global optimality of known locally convergent methods.
\end{abstract}

\section{INTRODUCTION}\label{sec:intro}
Real-world dynamical systems are often high-dimensional, resulting in significant computational and hardware demands making simulation and real-time implementation impractical. To address this challenge, model reduction has become an essential tool for simplifying system representations while preserving their key dynamic behavior. It has been extensively studied (e.g.,~\cite{study1, Antoulas}) and applied in areas such as control theory, fluid dynamics~\cite{challenge1}, and electrical engineering. In this note, the dynamics of a stable $n$th order SISO continuous-time LTI system $G(s)$ is described by:
\begin{equation}\label{sys1}
\left\{\begin{array}{l}
\dot{{x}}(t)={A} {x}(t)+{B} u(t), \\
y(t)=C {x}(t),
\end{array}\right.
\end{equation}
where ${A} \in \mathbb{R}^{n \times n},~{B} \in \mathbb{R}^{n \times 1}$ and ${C}\in \mathbb{R}^{1 \times n}$ and where ${x}(t) \in \mathbb{R}^n, ~u(t),~y(t) \in \mathbb{R}$ denote the state, input and output of the system. The transfer function is given by
\begin{equation}
G(s)=C(s {I}-{A})^{-1} {B}.
\end{equation}

A widely used technique, as discussed in \cite{Gugercin2008, Sato2018}, formulates the model reduction problem as an optimization task: given a stable system $G(s)$, the goal is to find a stable reduced-order model $G_m(s)$ that minimizes the $\mathcal{H}_2$ error norm \cite{Ahmad20102OM}:
\begin{equation}\label{H2NORM}
G_m=\underset{\substack{\hat{G}_m \text { is stable } \\ \operatorname{dim}\left(\hat{G}_m\right)\ =\ m}}{\arg } \min \left\|G-\hat{G}_m\right\|_{\mathcal{H}_2}
\end{equation}
where $\operatorname{dim}(\cdot)$ denotes the order of the transfer function. 

An effective method, Iterative Rational Krylov Algorithm (IRKA)  \cite{Gugercin2008}, provides an iterative procedure to compute locally optimal shifts for constructing the reduced-order system. This projection-based numerical method has been used across various applications (see, e.g., \cite{application3, application4}). To address the limitation of local optimality, the zero-solving approach proposed in \cite{zero-solvig2} generates interpolation points based on system zeros, aiming to account for all possible local optima. Alternative approaches, including balanced truncation \cite{balncedtruncation, balncedtruncation2}, Lanczos \cite{Lanczos1950AnIM,Lanczos}, and Arnoldi-based methods \cite{Arnoldi1951ThePO,Arnoldi,Arnoldi2}, provide solutions for large-scale model reduction. While these techniques are computationally efficient, they often struggle to ensure both high accuracy and stability in the reduced models. This trade-off highlights the ongoing need for robust and reliable model reduction methods tailored to specific application requirements. Despite several efforts in \cite{Antoulas2015, SPANOS1992897, Fulcheri1998}, the global optimality problem remains open.

In this paper, we focus on SISO systems and present a novel approach to determine a set of globally optimal interpolation points that minimize the $\mathcal{H}_2$ error in (\ref{H2NORM}) for $m=1,2$. Building on \cite{Ahmad20102OM}, we interpret the problem as one of maximizing the $\mathcal{H}_2$ norm within the set of all reduced order models satisfying the necessary conditions for optimality. Unlike iterative approaches and the method in \cite{Alsubaie}, which reformulates the problem into a multi-parameter matrix pencil problem, we use a semidefinite relaxation (SDR) approach and introduce an upper bound $\gamma$ to recast the problem into a convex optimization framework. Specifically, under necessary conditions and stability constraints imposed by the Routh–Hurwitz stability criterion, we explicitly express $\|G_m\|_{\mathcal{H}_2}^2$ as a function of interpolation points with an upper bound defined by $\gamma^2$.  The resulting SDR  method directly yields optimal interpolation points, thus overcoming convergence limitations of IRKA and significantly enhancing the applicability of existing iterative or Krylov-based methods.

The paper is organized as follows. Section~\ref{sec:Preliminaries} provides a brief review of the model reduction problem. The main results are provided and analyzed in Section \ref{sec:Approximation}. We propose a general approach to interpolation-based approximation and derive explicit results for $m=1$ and $m=2$. Numerical examples are detailed in Section \ref{sec:examples} and concluding remarks are given in Section \ref{sec:conclusion}. Finally, the appendix gives a proof of the optimality of our solution when $m=1$.

\textit{Notation.} $\mathbb{R}$ (resp., $\mathbb{C}$) denotes the set of real (complex) numbers, and $\mathbb{R}^{n\times m}$ (resp., $\mathbb{C}^{n\times m}$) the set of $n\times m$ real (complex) matrices.  For $X \in \mathbb{C}^{n\times m}$, its transpose, complex conjugate and conjugate transpose are denoted as $X^T$, $\bar{X}$ and $X^H = \bar{X}^T$, respectively. If $A=A^H$ then $A \succeq 0$ indicates that $A$ is positive semi-definite. The spectrum of a matrix $A\in\mathbb{C}^{n\times n}$ is denoted as $\lambda(A)$. The $n\times n$ identity matrix is denoted by $I_n$, with the subscript normally dropped if it can be inferred from the context. Given a stable system $G(s)$, $\|G\|_{\mathcal{H}_2}$ denotes the $\mathcal{H}_2$ norm. Additionally, for brevity, we define $\eta(x) := x + \bar{x}$.

\section{preliminaries}\label{sec:Preliminaries}
Consider the system (\ref{sys1}) as detailed in Section \ref{sec:intro}. Model reduction aims to construct a stable approximating system $G_{m}(s)$:  
\begin{equation}\label{sys2}
\left\{\begin{array}{l}
\dot{{x}}_m(t)={A}_m {x}_m(t)+{B}_m u(t), \\
y_m(t)={C}_m {x}_m(t),
\end{array}\right.
\end{equation}
with a reduced order $m$ ($m < n$), where ${A_m} \in \mathbb{R}^{m \times m}$ and ${B_m} \in \mathbb{R}^{m\times 1}$, ${C_m}\in \mathbb{R}^{1\times m}$. The transfer function of (\ref{sys2}) is 
\begin{equation}
G_m(s)=C_m(s {I_m}-{A_m})^{-1} {B_m}.
\end{equation}
\begin{definition}\cite{Gugercin2008}\label{def:h2}
A stable SISO system with a strictly proper transfer function $G(s)=C(s {I}-{A})^{-1} {B}$ has an expression of $\|G\|_{\mathcal{H}_2}^2$ such that
\begin{equation}\label{eq:Gh22}
\|G\|_{\mathcal{H}_2}^2  =\frac{1}{2 \pi} \int_{-\infty}^{+\infty} G^H(i \omega) G(i \omega) d \omega
 =C {P} {C^H},
\end{equation}
where ${P=P^H\succeq0}$ is the controllability Gramian and the unique solution to the Lyapunov equation: 
\begin{equation}\label{eq:Lyapunov}
{A}{P}+{P}{A}^H+{B} {B}^H=0.
\end{equation}
\end{definition}
The traditional model reduction strategy (\ref{H2NORM}) allows for multiple local solutions.
\begin{definition}\cite{Gugercin2008}
A stable order $m$ system $G_m$ is a local minimizer of the problem (\ref{H2NORM}) if it satisfies:
\begin{equation}\label{eq:cond1}
\left\|G-G_m\right\|_{\mathcal{H}_2} \leq\left\|G-\tilde{G}_m^{(\epsilon)}\right\|_{\mathcal{H}_2},
\end{equation}
for all feasible $\tilde{G}_m^{(\epsilon)}$ with $dim(\tilde{G}_m^{(\epsilon)})=m$ satisfying $\left\|G_m - \tilde{G}_m^{(\epsilon)}\right\|_{\mathcal{H}_2} \leq C\epsilon$, where $\epsilon > 0$ is sufficiently small and $C$ is a constant.

\end{definition}
\begin{definition}\cite{moment_matching}\label{def:moment}
The $j$-moment of system (\ref{sys1}) at $s \in \mathbb{C}$ is the complex number
\begin{equation}
\zeta_j(s)={C}\left(s{I}-{A}\right)^{-(j+1)} {B}.
\end{equation}
\end{definition} 
\vspace{1mm}
The point $s$ at which $\zeta_j(s)$ is computed is often referred to as an interpolation point. For our purposes, extending the equality to the first two moments is sufficient.
\begin{lem}[necessary condition]\cite{Gugercin2008}\label{lemma: necessary condtion}
For the stable SISO system $G(s)$ in (\ref{sys1}) of dimension $n$ let $G_m(s)={C_m}(s {I_m}-{A_m})^{-1} {B_m}$ be a local minimizer and assume that it has simple poles at $\widetilde{\lambda}_i$, where $i=1,2, \ldots, m$. Then, $G_m(s)$ interpolates both $G(s)$ and its first derivative at $-\widetilde{\lambda}_i$, $i=1, \dots m$:
\begin{equation}\label{eq:G=G'}
G_m\left(s_i\right)=G\left(s_i\right), \quad G_m^{\prime}\left(s_i\right)=G^{\prime}\left(s_i\right),
\end{equation}
where $s_i=-\widetilde{\lambda}_i$ for $i=1,2, \ldots, m$.
\end{lem} 

The set $S_m = \left\{s_1, s_2, \ldots, s_m\right\}$ corresponds to the negative eigenvalues of $A_m$, as stated above. A real realization of $G_m(s)$ is presented next, satisfying the moment matching conditions in (\ref{eq:G=G'}) for a given set $S_m$.
\begin{lem}
\label{def:Vm Wm}\cite{Alsubaie}
Let $s_i \in \mathbb{C}$ be distinct such that the set $S_m = \left\{s_1, s_2, \ldots, s_m\right\}$ is closed under complex conjugation and $S_m\cap \lambda(A) = \emptyset$. Define 
\begin{equation}\label{eq:V & W}
\begin{aligned}
V_m\left(S_m\right)&= \left(\prod_{j=1}^m\left(s_j I-A\right)^{-1}\right)\left[\begin{array}{llll}
B & A B & \cdots & A^{m-1} B
\end{array}\right]\\
W_m\left(S_m\right) & =\left(\left[\begin{array}{c}
C \\
C A \\
\vdots \\
C A^{m-1}
\end{array}\right]\left(\prod_{j=1}^m\left(s_j I-A\right)^{-1}\right)\right)^T
\end{aligned}
\end{equation}   
and assume that $T_m\left(S_m\right):=W_m\left(S_m\right)^T V_m\left(S_m\right)$ is non-singular. Then, 
\begin{eqnarray}\label{eq:Am}
A_m\left(S_m\right)&=&T_m\left(S_m\right)^{-1} W_m\left(S_m\right)^T A V_m\left(S_m\right)\label{eq:Bm}\\
B_m\left(S_m\right)&=&T_m\left(S_m\right)^{-1} W_m\left(S_m\right)^TB\label{eq:Cm}\\
C_m\left(S_m\right)&=&CV_m\left(S_m\right)
\end{eqnarray}
gives a real realization of $G_m(S_m)$ satisfying the interpolation condition in (\ref{eq:G=G'}).
\end{lem}
\begin{proof}
By Definition \ref{def:Vm Wm}, the rational interpolation is obtained through the Krylov projection method \cite{RUHE1994283, Grimme1997PhD, Lanczos}. Equation (\ref{eq:Am}) follows directly from computations in \cite[equ.~(3.7) and (3.8)]{Alsubaie}, while the necessary condition for $G_m(S_m)$ is established in \cite[Lem 3.2]{Gugercin2008}.
\end{proof}
\begin{rem}
Any $S_m$ that satisfies the necessary conditions in Lemma \ref{lemma: necessary condtion} is a fixed point of the transformation: $S_m = \lambda\left(-A_m(S_m)\right)$, where $A_m(S_m)$ is defined in \eqref{eq:Am}.
\end{rem}

However, selecting $S_m$ for constructing a local minimizer $G_m(s)$ presents a key challenge: ensuring that $\lambda(A_m) = -S_m$.

The following lemma, summarized from \cite{Ahmad20102OM}, provides the theoretical basis for the main contribution of this work.
\begin{lem}\cite{Ahmad20102OM}\label{lem:maxGm}
Let  $S_m=\left\{s_1, \ldots, s_m\right\} \subset \mathbb{C}^{+}$ be a fixed point of 
\begin{equation}\label{eq:fixed point}
S_m = \lambda\left(-A_m(S_m)\right),
\end{equation}
where $A_m(S_m)$ is defined in \eqref{eq:Am} and let $A_m,~B_m$ and $C_m$ be as defined in \eqref{eq:Am}-\eqref{eq:Cm}, respectively, with the dependence on $S_m$ dropped for convenience. Let $G_m(s)=C_m(sI_m-A_m)^{-1}B_m$,
so that $G_m(s)$ satisfies the necessary condition in Lemma~\ref{lemma: necessary condtion}. Then the following equation holds:
\begin{equation}\label{eq:global h2}
\left\|G-G_m\right\|_{\mathcal{H}_2}^2=\|G\|_{\mathcal{H}_2}^2-\left\|G_m\right\|_{\mathcal{H}_2}^2.
\end{equation}
\end{lem}
\vspace{2mm}
\begin{proof}
The proof follows from a derivation by using (\ref{eq:Gh22}) and the rational Arnoldi equations in \cite{FRANGOS2008342}.
\end{proof}

\section{Interpolated Approximation of Reduced Order System} \label{sec:Approximation}
In this section, we reformulate the $\mathcal{H}_2$ model reduction problem to reveal its underlying maximization form. This formulation forms the basis for deriving explicit approximations for $m = 1$ and $m = 2$ in the subsequent analysis.

\begin{thm}\label{thm:Gmmax}
Given a stable system $G(s)$ as defined in (\ref{sys1}), let $\hat{\mathcal{G}}_m$ denote the set of all transfer functions that satisfy the necessary conditions stated in Lemma~\ref{lemma: necessary condtion}. Then 
\begin{equation}\label{eq:ref_Gm}
\underset{\substack{G_m \text { is stable } \\ \operatorname{dim}\left(G_m\right)=m}}{\arg }\!\!\!\!\!\!\!\! \min \left\|G-G_m\right\|_{\mathcal{H}_2}\!=\!\!\!\!\!\underset{\substack{G_m \text { is stable } \\ \operatorname{dim}\left(G_m\right)=m
\\G_m\in\hat{\cal{G}}_m}}{\arg }\!\!\!\!\!\!\! \max \left\|G_m\right\|_{\mathcal{H}_2}.
\end{equation}
\end{thm}
\vspace{1mm}
\begin{proof}
This follows from Lemma \ref{lem:maxGm}, which implies that the system in $\hat{\cal{G}}_m$ with the maximum $\mathcal{H}_2$ norm is the solution to the minimum error problem as in (\ref{H2NORM}).
\end{proof}

Thus, solving (\ref{H2NORM}) is equivalent to identifying the system in $\hat{\mathcal{G}}_m$ with the maximum $\mathcal{H}_2$ norm.


To simplify the presentation, we begin by transforming the realization of $G_m(s)$ into a special diagonal form that preserves its dynamical behavior. The resulting approximate system is given by
\begin{equation}\label{eq:transfer}
\left[\begin{array}{c|c}
A_m & B_m \\\hline C_m &0\end{array}\right]
=
\left[
\begin{array}{ccc|c}
a_1 & \cdots & 0 & b_1 \\
0 & \ddots & 0 & \vdots \\
0 & \cdots & a_m & b_m \\ 
\hline
1 & \cdots & 1 & 0
\end{array}
\right].
\end{equation}
Note that $\lambda(A_m)=\{a_1,\ldots,a_m\}$. Note also that there is no loss of generality in assuming this form since Lemma~\ref{lemma: necessary condtion} assumes that $G_m(s)$ has simple poles.

For convenience, we also 
characterize the set of interpolating points $S_m$ through the following change of variables:
\begin{equation}\label{eq:si2pi}
\begin{aligned}
p_0 & =1 \\
p_1 & =\sum_{1 \leq j \leq m} s_j \\
& \vdots \\
p_k & =\sum_{1 \leq j_1<\cdots<j_k \leq m} s_{j_1} \cdots s_{j_k} \\
& \vdots \\
p_m & =s_1 \cdots s_m
\end{aligned}
\end{equation}
Note that the $p_i$'s are real and that the $s_i$'s can be obtained from $p_i$'s as the solution of the polynomial equation:
\begin{equation}\label{eq:poly equation}
s^m - p_1 s^{m-1} + \cdots + (-1)^{m-1} p_{m-1}s + (-1)^m p_m = 0.
\end{equation}
\subsection{First order approximation}
Consider the model reduction problem when $m=1$. Let $G_1(s)=C_1(sI_1-A_1)^{-1}B_1\in\hat{\mathcal{G}}_1$, with the form given in \eqref{eq:transfer}, be a local minimizer for $G(s)$ and let $s_1$ be the interpolation point. Then $s_1$ is a fixed point of \eqref{eq:fixed point} if $a_1=-s_1$. Since $G_1(s_1)=G(s_1)$, we have
\vspace{-2mm}
\begin{equation}\label{eq:sysm=1}
C_1\left(s_1 I_1-A_1\right)^{-1} B_1=\frac{b_1}{2 s_1}=C\left(s_1 I_n-A\right)^{-1} B.
\end{equation}
Since $G_1(s)$ is a scalar transfer function, we can assume that $a_1, b_1 \in \mathbb{R}$. 

Next, we derive an explicit expression for $\|G_1\|_{\mathcal{H}_2}^2$ via the controllability Gramian using Definition~\ref{def:h2}.
A unique solution to (\ref{eq:Lyapunov}):
\begin{equation}
2 a_1 P_1+b_1^2=0,    
\end{equation}
is given by
$
P_1=-\frac{b_1^2}{2 a_1}=\frac{b_1^2}{2 s_1}.
$ 
It follows from \eqref{eq:Gh22} and \eqref{eq:sysm=1} that
\begin{equation}
\begin{aligned}
\| G_1 \|_{\mathcal{H}_2}^2&=C_1 P_1 C_1^T=\frac{b_1^2}{2 s_1}\\
&=2 C\left(s_1 I-A\right)^{-1} s_1 B B^T\left({s}_1 I-A\right)^{-T} C^T.
\end{aligned}
\end{equation}
An alternative representation of $\|G_1\|_{\mathcal{H}_2}^2$ using the change of variables in \eqref{eq:si2pi} gives
\begin{equation}\label{eq:f(p)m=1 prior}
\hspace{-0.21cm}
\|G_1\|_{\mathcal{H}_2}^2\!=\!2 C\!\left(p_1 I\!-\!A\right)^{-1}\!\!\!\left(p_1 B B^T\!\right)\!\!\left({p_1} I\!-\!A\right)^{-T}\!\! C^T\!\!=:\!\!f(p_1),
\end{equation}
where  $p_1=s_1$. We can rewrite $f(p_1)$ as
$$
\begin{aligned}
f(p_1) =&2 C(\overbrace{p_1 I-A}^{\mathcal{A}_1})^{-1}\left(p_1 B B^T\right) \mathcal{A}_1^{-T} C^T \\
=&2 C(I-\overbrace{p_1 I}^{\Delta} \overbrace{A^{-1}}^{T_4})^{-1}{p_1}(\overbrace{A^{-1} B B^T A^{-T}}^{T_6}) \mathcal{A}_1^{-T} C^T \\
=&2 C\left(I-\Delta T_4\right)^{-1}\left(\Delta T_6\right)\left(I-\Delta T_4\right)^{-T} C^T \\
=&2 C \Delta\left(I-T_4 \Delta\right)^{-1} T_6\left(I+T_4^T\left(I-T_4 \Delta\right)^{-T} \Delta^T\right)^T C^T \\
=&C \overbrace{\Delta\left(I-T_4 \Delta\right)^{-1}}^{\widehat{\Delta}} T_6 C^T+C T_6^T\left(I-T_4 \Delta\right)^{-T} \Delta^T C^T \\
&+C \Delta\left(I-T_4 \Delta\right)^{-1} \overbrace{\left(T_6 T_4^T+T_4 T_6^T\right)}^{-T_5}\left(I-T_4 \Delta\right)^{-T} \Delta^T C^T \\
=&C \widehat{\Delta} T_6 C^T+C T_6^T \widehat{\Delta}^T C^T-C \widehat{\Delta} T_5 \widehat{\Delta}^T C^T.
\end{aligned}
$$
Let $\gamma^2$ be an upper bound on $f(p_1)$ for all $p_1>0$. Then
\begin{equation}\label{eq:001}
\gamma^2-f(p_1) =\gamma^2+T_2 \widehat{\Delta} T_3+T_3^T \widehat{\Delta}^T T_2^T+T_2 \widehat{\Delta} T_5 \widehat{\Delta}^T T_2^T,
\end{equation}
where $T_2=-C, T_3=T_6 C^T$.

The stability of $G_1(s)$ imposes the constraint on the interpolation point that $s_1>0$ (equivalently, $p_1>0$), which the next proposition incorporates, providing a convex semi-definite optimization.
\begin{prop}\label{prop:m=1}
Consider the model reduction problem for $m = 1$ and let $\gamma^2$ be an upper bound on $f(p_1)$ for all $p_1>0$. Then, $\gamma^2$ can be minimized subject to the LMIs \begin{equation}\label{C_1}
S + S^T \succeq 0,\quad L \succeq 0,
\end{equation}
where $S \in \mathbb{R}^{n \times n}$ and $L=L^T\in \mathbb{R}^{(n+1)\times(n+1)}$ is defined as
\[L = \begin{bmatrix}
\gamma^2 & T_3^T - T_2 S^T \\
T_3 - S T_2^T & -S T_4^T - T_4 S^T + T_5
\end{bmatrix}.\]
\end{prop}
\vspace{1.5mm}
\begin{proof}
Note that the reduced order model is feasible if and only if $s_1 > 0$, or equivalently, $p_1 > 0$. To enforce this constraint, let $S \in \mathrm{R}^{n \times n}$ be a slack variable with $S+S^T \succeq 0$, which yields
\begin{equation}\label{eq:S+S'}
p_1\left(S+S^T\right) \succeq 0,~\forall p_1>0.
\end{equation} 
Since $\Delta\!=\!p_1 I$ and $\Delta\!=\!\left(I\!+\!\hat{\Delta} T_4\right)^{-1} \hat{\Delta}$, this inequality is equivalent to the following inequalities satisfied for all $p_1\!>\!0$:
\begin{gather}\label{eq:reformulated S+S'}
    \Delta S + S^T \Delta^T \succeq 0, \notag \\
    \left(I+\hat{\Delta} T_4\right)^{-1} \hat{\Delta} S + S^T \hat{\Delta}^T\left(I+T_4^T \hat{\Delta}^T\right)^{-1} \succeq0, \notag \\
    \hat{\Delta} S + S^T \hat{\Delta}^T + \hat{\Delta}\left(S T_4^T + T_4 S^T\right) \hat{\Delta}^T \succeq 0, \notag \\
    \label{eq:lasteq}
    T_2 \hat{\Delta} S T_2^T + T_2 S^T \hat{\Delta}^T T_2^T + T_2 \hat{\Delta}\left(S T_4^T + T_4 S^T\right) \hat{\Delta}^T T_2^T \succeq 0.
\end{gather}
A simple manipulation shows that
$$
\gamma^2 \! -\! f(p_1)\! =\! \eta(w)\!+\!T_2 \hat{\Delta}\left(S T_4^T\!+\!T_4 S^T\right) \hat{\Delta}^T T_2^T\! +v L v^T
$$
where $w=T_2 \hat{\Delta} S T_2^T$,
$v = \left[\begin{array}{cc}
1 & T_2 \hat{\Delta}
\end{array}\right]$. It follows from (\ref{eq:lasteq}) and (\ref{eq:001}) that 
\vspace{-2mm}
\begin{equation}\label{eq:002}
\gamma^2 \! -\! f(p_1)\! =\! \overbrace{\eta(w)\!+\!T_2 \hat{\Delta}\left(S T_4^T\!+\!T_4 S^T\right) \hat{\Delta}^T T_2^T}^{\succeq 0}\! +v L v^T\!\! \geq 0,
\end{equation}
 for all $p_1\!>\!0$ when $L\!\succeq\!0$. In turn, the upper bound of $\|G_1\|_{\mathcal{H}_2}^2$ is obtained by minimizing $\gamma^2$ under the LMIs \eqref{C_1}. 
\end{proof}


\subsection{Second order approximation}
Consider $G_m(s)$ with $m\!=\!2$, where the interpolating points are $s_1, s_2\! \in\! \mathbb{R}$ for the real case and $s_1, s_2\! \in \!\mathbb{C}$, $s_2\!=\!\bar{s}_1$ for the complex case. Since $s_1,s_2$ are interpolating points, then
\begin{equation}\label{eq:sysm=2}
\left[\begin{array}{c}
C_m\left(s_1 I-A_m\right)^{-1} B_m \\
C_m\left(s_2 I-A_m\right)^{-1} B_m
\end{array}\!\!\right]\!\!=\!\!\left[\begin{array}{l}
C\left(s_1 I-A\right)^{-1} B \\
C\left(s_2 I-A\right)^{-1} B
\end{array}\!\!\right]\!\!:=\!\!\left[\!\!\begin{array}{l}
g_1 \\
g_2
\end{array}\!\!\right]\!\! .
\end{equation}
Since $S_m=-\lambda\left(A_m\right)$, $s_1=-a_1$, $s_2=-a_2$. Hence
\begin{equation}\label{eq:sysm=2.1}
\left[\begin{array}{cc}
\frac{1}{2 s_1} & \frac{1}{s_1+s_2} \\
\frac{1}{s_1+s_2} & \frac{1}{2 s_2}
\end{array}\right]\left[\begin{array}{l}
b_1 \\
b_2
\end{array}\right]=\left[\begin{array}{l}
g_1 \\
g_2
\end{array}\right].
\end{equation}
Next, we evaulate the controllability Gramian explicitly.
\begin{lem}\label{lemma:pm m=2}
For $m=2$, consider $G_m(s)$ in the form of (\ref{eq:transfer}). Define 
$$P_m=P_m^H=\begin{bmatrix}
P_1 & P \\ \bar{P} & P_2
\end{bmatrix},$$
satisfying (\ref{eq:Gh22}). Then
\begin{equation}\label{eq:pmm=2}
P_m=\left[\begin{array}{cc}
\frac{b_1 \bar{b}_1}{s_1+\bar{s}_1} & \frac{b_1 \bar{b}_2}{s_1+\bar{s}_2} \\
\frac{\bar{b}_1 b_2}{s_2+\bar{s}_1} & \frac{b_2 \bar{b}_2}{s_2+\bar{s}_2}
\end{array}\right].
\end{equation} 
\end{lem}
\vspace{1.5mm}
\begin{proof}
The computation of $P_m$ follows straightforward mathematical steps. For $m=2$, (\ref{eq:Gh22}) gives that
\begin{equation}
\operatorname{diag}(a_1, a_2)  
\!\begin{bmatrix}
P_1 & \!\!P \\ \bar{P} & \!\!P_2
\end{bmatrix}\!
+\!
\begin{bmatrix}
P_1 & \!\!P \\ \bar{P} & \!\!P_2
\end{bmatrix}
\operatorname{diag}(\bar{a}_1, \bar{a}_2)
\!+\!
\mathbf{b} \mathbf{b}^H
\!=\! 0,
\end{equation}
where $\mathbf{b} = \begin{bmatrix} b_1 & b_2 \end{bmatrix}^T$, $a_i$ and its conjugates are replaced by $s_i$ and $\bar{s}_i$ respectively, yielding:
\begin{equation}
\left[\begin{array}{cc}
\left(s_1+\bar{s}_1\right) P_1 & \left(s_1+\bar{s}_2\right) P \\
\left(s_2+\bar{s}_1\right) \bar{P} & \left(s_2+\bar{s}_2\right) P_2
\end{array}\right]=\left[\begin{array}{cc}
b_1 \bar{b}_1 & b_1 \bar{b}_2 \\
\bar{b}_1 b_2 & b_2 \bar{b}_2
\end{array}\right].
\end{equation}
Thus, the matrix $P_m$ is directly obtained as (\ref{eq:pmm=2}).
\end{proof}

It follows from \eqref{eq:Gh22} that
\begin{equation}\label{eq:Gmm=2}
\left\|G_m\right\|_{\mathcal{H}_2}^2\!=\!\frac{b_1 \bar{b}_1}{\left(s_1\!+\!\bar{s}_1\right)}\!+\!\frac{b_1 \bar{b}_2}{\left(s_1\!+\!\bar{s}_2\right)}\!+\!\frac{\bar{b}_1 b_2}{\left(s_2\!+\!\bar{s}_1\right)}\!+\!\frac{b_2 \bar{b}_2}{\left(s_2\!+\!\bar{s}_2\right)}\cdot
\end{equation}
Next, we evaluate $\left\|G_m\right\|_{\mathcal{H}_2}^2$ explicitly.
\begin{lem}\label{lemma:f(p)m=2}
Denote $\left\|G_m\right\|_{\mathcal{H}_2}^2$ as $f(s_1,s_2)$. Then
\vspace{-0.5mm}
\begin{equation}\label{eq:f(p)m=2}
f(s_1,s_2)\!=\!2(s_1\!+\!s_2) C \mathcal{A}_2^{-1}\left(s_1s_2 \mathbf{B}\!+\! A \mathbf{B} A^T\right) \mathcal{A}_2^{-T} C^T,
\vspace{-0.5mm}
\end{equation}
where $\mathcal{A}_2=\left(s_1 I-A\right)\left(s_2 I-A\right),\mathbf{B}=B B^T$.
\end{lem}
\begin{proof}
Note first that since $G_m(s)$ has the form in \eqref{eq:transfer}, we can assume that $a_1\ne a_2$ (otherwise, $G_m(s)$ is a first order system), so $s_1\ne s_2$. It follows from (\ref{eq:sysm=2.1}) that
\begin{align*}
\begin{bmatrix}
b_1 \\
b_2
\end{bmatrix}=\frac{4 s_1 s_2\left(s_1+s_2\right)^2}{\left(s_1-s_2\right)^2}\left[\begin{array}{l}
\frac{g_1}{2 s_2}\!-\frac{g_2}{s_1+s_2} \\
\frac{g_2}{2 s_1}\!-\frac{g_1}{s_1+s_2}
\end{array}\right].
\end{align*}
Let $M=\frac{4 s_1 s_2\left(s_1+s_2\right)^2}{\left(s_1-s_2\right)^2}$ and define $\psi(x)=x\bar{x},$ for any $ x \in \mathbb{C}$. A direct computation yields that
$$
\begin{aligned}
&b_1 \bar{b}_1 = M^2 \left( 
    \frac{g_1 \bar{g}_1}{4 s_2 \bar{s}_2} 
    + \frac{\bar{g}_2 g_2}{\psi(s_1+s_2)}
    - \eta \left( \frac{g_2 \bar{g}_1}{2 \bar{s}_2 (s_1+s_2)} \right)
\right), \\
& b_2 \bar{b}_2=M^2\left(\frac{g_1 \bar{g}_1}{\psi(s_1+s_2)}+\frac{g_2 \bar{g}_2}{4 \bar{s}_1 s_1}-\eta\left(\frac{g_2 \bar{g}_1}{2 s_1\left(\bar{s}_1+\bar{s}_2\right)}\right)\right), \\
& b_1 \bar{b}_2 = M^2X_{b_1b_2},\\
&X_{b_1b_2}=\frac{g_1 \bar{g}_2}{4 s_2 \bar{s}_1} + \frac{g_2 \bar{g}_1}{\psi(s_1+s_2)}- \frac{g_1 \bar{g}_1}{2 s_2 (\bar{s}_1 + \bar{s}_2)}- \frac{g_2 \bar{g}_2}{2 \bar{s}_1 (s_1 + s_2)}.
\end{aligned}
$$
Note that if $s_1$, $s_2\in \mathbb{R}$, $s_1=\bar{s}_1, s_2=\bar{s}_2$ and if $s_1$, $s_2\in \mathbb{C}$, $\bar{s}_2=s_1, \bar{s}_1=s_2$. By explicit expansion of (\ref{eq:Gmm=2}), we get
\begin{equation}\label{real}
f(s_1,s_2)\!=\!\frac{2\left(s_1\!+\!s_2\right)(\left(s_1 g_1\!-\!s_2 g_2\right)^2\!+\!s_1 s_2\left(g_1\!-\!g_2\right)^2)}{\left(s_2\!-\!s_1\right)^2},
\end{equation}
for the real case, and 
\vspace{-0.5mm}
\begin{equation}\label{complex}
\hspace{-4mm}f(s_1,s_2)\!=\!\frac{-2\left(s_1\!+\!s_2\right)\left(\psi\left(s_1 g_1\!-\!s_2 g_2\right)+s_1 s_2\psi\left(g_1\!-\!g_2\right)\right)}{\left(s_2\!-\!s_1\right)^2}
\end{equation}
for the complex case. Reconstruct $g_1$ and $g_2$,
$$
\begin{aligned}
& g_1=C\left[\left(s_1I-A\right)^{-1}\left(s_2I-A\right)^{-1}\right]\left(s_2I-A\right) B, \\
& g_2=C\left(s_1I-A\right)\left[\left(s_1I-A\right)^{-1}\left(s_2I-A\right)^{-1}\right]B.
\end{aligned}
$$
Substituting these expressions into (\ref{real}) and (\ref{complex}) yields the same result as in (\ref{eq:f(p)m=2}). Therefore, (\ref{eq:f(p)m=2}) is proven.
\end{proof}

From \eqref{eq:si2pi}, we obtain $p_1 = s_1 + s_2$ and $p_2 = s_1 s_2$. An alternative expression for  $\|G_m\|_{\mathcal{H}_2}^2$ is given by
\begin{equation}\label{eq:f(p)m=2 prior}
f(p)=2 C \mathcal{A}_2^{-1}\left(p_1 p_2 B B^T+p_1 A B B^T A^T\right) \mathcal{A}_2^{-T} C^T,
\end{equation}
where $\mathcal{A}_2\!=\!p_2 I\!-\!p_1 A\!+\!A^2$ and  $p$ denotes the pair $(p_1,p_2)$. Rewriting (\ref{eq:f(p)m=2 prior}) gives 
$$
\begin{aligned}
f(p)=& 2 C(\overbrace{p_2 I-p_1 A+A^2}^{\mathcal{A}_2})^{-1}\left(p_1 p_2 \mathbf{B}+p_1 A \mathbf{B} A^T\right) \mathcal{A}_2^{-T} C^T \\
=& 2 C(I-\overbrace{\left[\begin{array}{ll}
\underbrace{p_1 I}_{\Delta_1 I} & \underbrace{p_2 I}_{\Delta_2}
\end{array}\right]}^{\Delta} \overbrace{\left[\begin{array}{c}
A^{-1} \\
-A^{-2}
\end{array}\right]}^{T_4})^{-1}\\
& \times (\Delta\overbrace{\left[\begin{array}{c}
A^{-1} B B^T A^{-T} \\
0
\end{array}\right]}^{T_6}\\
&+\Delta\overbrace{\left[\begin{array}{cc}
0 & A^{-2} B B^T A^{-2 T} \\
0 & 0
\end{array}\right]}^{T_7} \Delta^T) \mathcal{A}_2^{-T} C^T\\
=& 2 C\left(I-\Delta T_4\right)^{-1}\left(\Delta T_6+\Delta T_7 \Delta^T\right)\left(I-\Delta T_4\right)^{-T} C^T\\
=& C \widehat{\Delta} T_6 C^T+C T_6^T \widehat{\Delta}^T C^T-C \widehat{\Delta} T_5 \widehat{\Delta}^T C^T,
\end{aligned}
$$
where $\mathbf{B}=B B^T,~T_5=-\left(T_6 T_4^T+T_4 T_6^T+T_7+T_7^T\right)$ and $\widehat{\Delta}=\left(I-\Delta T_4\right)^{-1}\Delta$. 

Let $\gamma^2$ denote an upper bound on $f(p)$ for all $p$. Then
\begin{equation}\label{eq:003}
\begin{aligned}
\gamma^2-f(p) =\gamma^2+T_2 \widehat{\Delta} T_3+T_3^T \widehat{\Delta}^T T_2^T+T_2 \widehat{\Delta} T_5 \widehat{\Delta}^T T_2^T,
\end{aligned}
\end{equation}
where $T_2=-C, T_3=T_6 C^T$.

The following similarly ensures the stability of $G_m(s)$ for $m=2$ and provides an SDR method to compute the minimum value of the upper bound $\gamma^2$.
\begin{prop}\label{prop:m=2}
Consider the model reduction problem for \(m = 2\). Let \(\gamma^2\) be an upper bound on $f(p)$ for $p=(p_1,p_2)$, $p_1>0,~p_2>0$. Then, \(\gamma^2\) can be minimized subject to 
\begin{equation}\label{C_2}
S_1+S_1^T \succeq0,~~S_2+S_2^T\succeq 0,~~G_{12}+G_{12}^T \succeq 0,~~L \succeq 0,
\end{equation}
where $S \in \mathbb{R}^{2n \times n}$, $G \in \mathbb{R}^{2n \times 2n}$ and \(L=L^T\in \mathbb{R}^{(2n+1)\times(2n+1)}\) are defined as
\[S=\left[\begin{array}{c}
S_1 \\
S_2
\end{array}\right], \quad G=\left[\begin{array}{cc}
0 & G_{12} \\
G_{12}^T & 0
\end{array}\right],\]
\[ L=\left[\begin{array}{cc}
\gamma^2 & T_3^T-T_2 S^T \\
T_3-S T_2^T & -S T_4^T-T_4 S^T+T_5-G
\end{array}\right].\]
\end{prop}
\vspace{2mm}
\begin{proof}
By the Routh-Hurwitz stability criterion, $p_1 > 0$ and $p_2 > 0$, derived from $s^2 + p_1 s + p_2 = 0$, where $s$ is replaced by $-s$ in (\ref{eq:poly equation}) due to $S_m = \lambda(-A_m(S_m))$. The redundant condition $p_1 p_2 > 0$ is further imposed.
Set $S_1+S_1^T \succeq 0$ and $S_2+S_2^T \succeq 0$, leading to 
$$
\left[\begin{array}{ll}
\Delta_1 & \Delta_2
\end{array}\right]\left[\begin{array}{c}
S_1 \\
S_2
\end{array}\right]+\left[\begin{array}{ll}
S_1^T & S_2^T
\end{array}\right]\left[\begin{array}{c}
\Delta_1^T \\
\Delta_2^T
\end{array}\right] \succeq 0
$$
for all $\Delta_1$ and $\Delta_2$.
Set $G_{12}+G_{12}^T \succeq 0$ to capture the quadratic condition $p_1p_2>0$, then
$$
\left[\begin{array}{ll}
\Delta_1 & \Delta_2
\end{array}\right]\left[\begin{array}{cc}
0 & G_{12} \\
G_{12}^T & 0
\end{array}\right]\left[\begin{array}{c}
\Delta_1^T \\
\Delta_2^T
\end{array}\right] \succeq 0
$$
for all $\Delta_1$ and $\Delta_2$. Combining both conditions, stability is captured by
\begin{equation}\label{eq:SG_m=2}
\Delta S + S^T \Delta^T + \Delta G \Delta^T \succeq 0
\end{equation}
for all $\Delta$.

Following (\ref{eq:reformulated S+S'}), (\ref{eq:003}) is obtained in the form of LMIs,
\begin{equation}\label{eq:004}
\begin{aligned}
\gamma^2-f(p)=\overbrace{\eta(w)+T_2 \hat{\Delta}\left(S T_4^T\!+T_4 S^T\!+G\right) \hat{\Delta}^T T_2^T}^{\succeq 0}\! +v L v^T\geq 0
\end{aligned}
\end{equation}
where 
$w\!=\!T_2 \hat{\Delta} S T_2^T\!$, $v \!=\!\! \left[\!\!\begin{array}{cc}
1 & \!\!\!T_2 \hat{\Delta}
\end{array}\!\!\right]$
and the first part is the congruence transformation of (\ref{eq:SG_m=2}). Thus the upper bound is obtained by minimizing $\gamma^2$  subject to the LMIs  \eqref{C_2}.
\end{proof}
\subsection{Optimal interpolation points}
The interpolating points can be further calculated once we obtain the optimal solution of the proposed SDR problem. To precisely determine the upper bound, we consider the ideal case where $\gamma^2-f(p)=0$, with $L \succeq0$ and where we assume that the zero eigenvalue of $L$ is simple. This leads to
\begin{equation}
L\left[\begin{array}{l}
1 \\
\mathrm{X}
\end{array}\right]=0,
\end{equation}
where $\left[1, X^T\right]^T$ is the normalized eigenvector of $L$ corresponding to the zero eigenvalue and $X=\hat{\Delta}^T T_2^T$. Using the definitions of $\hat{\Delta}$ and $T_2$, $X=\Delta^T\left(T_4^T X - C^T\right).$
Set $Z=T_4^T X-C^T$. Then all optimal interpolation points can be obtained by:
\begin{equation}\label{eq:solutionDelta}
\Delta^T=X . / Z 
\end{equation}
where $./$ denotes elementwise division. Thus, $p_1$ and $p_2$ (hence $s_1$ and $s_2$ via the roots of \eqref{eq:poly equation}) are directly obtained from the entries. In summary, the proposed approach is presented in Algorithm~\ref{alg:interpolation_points}.
\begin{algorithm}[ht]
\caption{$\mathcal{H}_2$ Optimal Interpolation Points and Reduced Order Model for First and Second Order Approximation of Linear Systems}
\label{alg:interpolation_points}
\textbf{Input:} $G(s) = C(sI - A)^{-1}B$, $n$, $m$\\
\textbf{Output:} Set of optimal interpolating points $S_m = \{s_1, \dots, s_m\}$ and optimal reduced order model $G_m(s)$\\
\textbf{Step 1: Problem Reformulation}\\
Compute the function $f(p)$ using (\ref{eq:f(p)m=1 prior}) if $m=1$ and (\ref{eq:f(p)m=2 prior}) if $m=2$. Further characterize the problem using (\ref{eq:001}) and (\ref{eq:003}), where $\gamma^2$ is an upper bound.\\
\textbf{Step 2: Solve the SDR Problem}\\
Determine the smallest upper bound $\gamma$ by solving the LMI optimization as described in Proposition \ref{prop:m=1} and Proposition \ref{prop:m=2}.\\
\textbf{Step 3: Computation of  Optimal Interpolation Points}\\
Compute $\Delta$ using equation (\ref{eq:solutionDelta}). If $m=1$, set $\Delta = p_1 I$ with $p_1 = s_1$, and assign $S_m = \{s_1\}$. If $m=2$, set $\Delta = [p_1 I, p_2 I]$ where $p_1 = s_1 + s_2$ and $p_2 = s_1 s_2$, then assign $S_m = \{s_1, s_2\}$ as the roots of the polynomial \eqref{eq:poly equation} .\\
\textbf{Step 4: Computation of  Optimal Reduced Order Model}\\
Use Lemma~\ref{def:Vm Wm} to obtain a real realization for $G_m(s)$.
\end{algorithm}
\section{NUMERICAL EXAMPLE}\label{sec:examples}
We demonstrate the global optimality of the proposed SDR method through the following four cases and compare the results with those obtained using IRKA~\cite{Gugercin2008}, where the example system is taken from~\cite{SPANOS1992897}, and the zero-solving method~\cite{zero-solvig2, Ahmad20102OM}. All cases are implemented in MATLAB using CVX solvers~\cite{grant2014cvx}, including SDPT3, MOSEK, and SeDuMi. Consider the following transfer functions: 
\begin{align*}
G_1(s)\!&=\!\frac{s^2+15 s+50}{s^4+5 s^3+33 s^2+79 s+50},\\
G_2(s)\!&=\!\frac{-1.986 s^2+19.17 s-0.1606}{s^3+4.857 s^2+14.08 s+23.02},\\
G_3(s)\!&=\!\frac{-1.3369 s^3-4.8341 s^2-47.5819s-42.7285}{s^4+17.0728 s^3+84.9908 s^2+122.4400s+59.9309},\\
G_4(s)\!&=\!\frac{-1.2805 s^3-6.2266 s^2-12.8095s-9.3373}{s^4+3.1855 s^3+8.9263 s^2+12.2936s+3.1987}.
\end{align*}

The reduced-order system approximations are constructed using fixed points from SDR. The selected optimal shifts are listed in Table~\ref{tab:G1},\ref{tab:G2},\ref{tab:G3},\ref{tab:G4} comparing all locally optimal points obtained using the zero-solving method and IRKA. For all these cases, the corresponding relative $\mathcal{H}_2$ errors
$$
\left\|G(s) - G_m(s)\right\|_{\mathcal{H}_2}/\|G(s)\|_{\mathcal{H}_2},
$$
are reported in Table~\ref{tab:h2_error_comparison} to assess the effectiveness of the algorithm in minimizing the $\mathcal{H}_2$ error.

The results clearly demonstrate that the proposed approach achieves globally optimal solutions, as the SDR method yields the same fixed points that serve as global minimizers, consistent with the findings in \cite{zero-solvig2}. For $m=2$, the deviation is less than 1\%, attributed to the precision of the CVX solver. When the optimal shifts from the SDR method are used to initialize IRKA$_1$, the final interpolation points converge exactly to the inputs. In contrast, IRKA$_2$, initialized with random shifts via the \texttt{randn} command, fails to reach the global optimum, despite some results appearing numerically close. Notably, for $G_3(s)$, IRKA$_2$ converges to a local optimum matching the zero-solving result. These findings support the claim that the proposed approach successfully identifies globally optimal interpolation points and demonstrates the global optimality of IRKA in these cases. More importantly, as evidenced by Table~\ref{tab:h2_error_comparison}, the proposed method yields the lowest error across all cases, confirming its numerical accuracy.
\begin{table}[htbp]
    \centering
    \caption{Comparison for $G_1(s)$.}
    \label{tab:G1}
    \begin{tabular}{|c|c|c|}
        \hline
        \textbf{Method type} & \textbf{Optimal Shifts ($m=1$)} & \textbf{Optimal Shifts ($m=2$)} \\
        \hline
        SDR&\{0.5762\}  &  \{4.1936, 1.1538\}\\
        \hline
        IRKA$_1$&\{0.5762\} & \{4.1936, 1.1538\} \\
        \hline
        IRKA$_2$& \{0.5780\}& \{4.1933, 1.1539\} \\
        \hline
        Zero-solving & \{0.5762\} & \{4.1935, 1.1539\}\\
        \hline
    \end{tabular}
\end{table}
\begin{table}[htbp]
    \centering
    \caption{Comparison for $G_2(s)$.}
    \label{tab:G2}
    \begin{tabular}{|c|c|c|}
        \hline
        \textbf{Method type} & \textbf{Optimal Shifts ($m=1$)} & \textbf{Optimal Shifts ($m=2$)} \\
        \hline
        SDR&\{2.1364\}  &  \{0.6935 $\pm$ 3.2772i\}\\
        \hline
        IRKA$_1$&\{2.1364\} & \{0.6935 $\pm$ 3.2772i\} \\
        \hline
        IRKA$_2$&\{0.0405\} &\{0.6934 $\pm$ 3.2771i\} \\
        \hline
                          & \{2.1364\}(global)& \\
       Zero-solving       &  \{0.0028\}(local) &\{0.6935 $\pm$ 3.2772i\}\\
                          &  \{36.2325\}(local)&\\
        \hline
    \end{tabular}
\end{table}
\begin{table}[ht!]
    \centering
    \caption{Comparison for $G_3(s)$.}
    \label{tab:G3}
    \begin{tabular}{|c|c|c|}
        \hline
        \textbf{Method type} & \textbf{Optimal Shifts ($m=1$)} & \textbf{Optimal Shifts ($m=2$)} \\
        \hline
        SDR&\{0.7007\}  & \{0.7051, 39.2818\}\\
        \hline
        IRKA$_1$&\{0.7007\} & \{0.7051, 39.2818\} \\
        \hline
        IRKA$_2$& \{0.6987\} & \{0.8261 $\pm$ 0.6577i\} \\
        \hline
    \multirow{2}{*}{\centering Zero-solving}&\multirow{2}{*}{\centering \{0.7007\}}&\{0.7051, 39.2800\}(global) \\
           &   &\{0.8261 $\pm$ 0.6577i\}(local)\\   
        \hline
    \end{tabular}
\end{table}
\begin{table}[ht!]
    \centering
    \caption{Comparison for $G_4(s)$.}
    \label{tab:G4}
    \begin{tabular}{|c|c|c|}
        \hline
        \textbf{Method type} & \textbf{Optimal Shifts ($m=1$)} & \textbf{Optimal Shifts ($m=2$)} \\
        \hline
        SDR&\{0.7828\}  & \{0.2030, 1.2052\}\\
        \hline
        IRKA$_1$&\{0.7828\} & \{0.2030, 1.2052\} \\
        \hline
        IRKA$_2$&\{0.7879\} & \{0.2076, 0.7527\}\\
        \hline
    \multirow{2}{*}{\centering Zero-solving}&\multirow{2}{*}{\centering \{0.7828\}}&\{0.2030, 1.2052\}(global) \\
           &   &\{6.3628, 1.1692\}(local)\\   
        \hline
    \end{tabular}
\end{table}
\begin{table}[ht!]
    \centering
    \caption{Comparison of the relative $\mathcal{H}_2$ errors.}
    \label{tab:h2_error_comparison}
    \begin{tabular}{|c|c|c|c|}
        \hline
        \textbf{Model} &m& \textbf{SDR} & \textbf{IRKA}  \\
        \hline
        $G_1(s)$&1& $4.8175\times 10^{-1}$ &$4.8175\times10^{-1}$ \\
        \hline
         $G_1(s)$&2& $2.4427\times10^{-1}$ &$2.4427\times10^{-1}$ \\
        \hline
         $G_2(s)$&1& $9.3389\times10^{-1}$ &$1.0000$ \\
        \hline
         $G_2(s)$&2& $4.3557\times10^{-1}$ &$4.3557\times10^{-1}$ \\
        \hline
         $G_3(s)$&1& $3.2235\times10^{-1}$ &$3.2235\times10^{-1}$ \\
        \hline
         $G_3(s)$&2& $2.6760\times10^{-1}$ &$2.9978\times10^{-1}$ \\
        \hline
         $G_4(s)$&1& $3.5992\times10^{-1}$ &$3.5993\times10^{-1}$ \\
        \hline
         $G_4(s)$&2& $3.2707\times10^{-1}$ &$3.6280\times10^{-1}$ \\
        \hline
    \end{tabular}
\end{table}

\begin{figure}[ht]
    \centering
    \begin{subfigure}{0.239\textwidth}
        \centering
        \includegraphics[width=\linewidth]{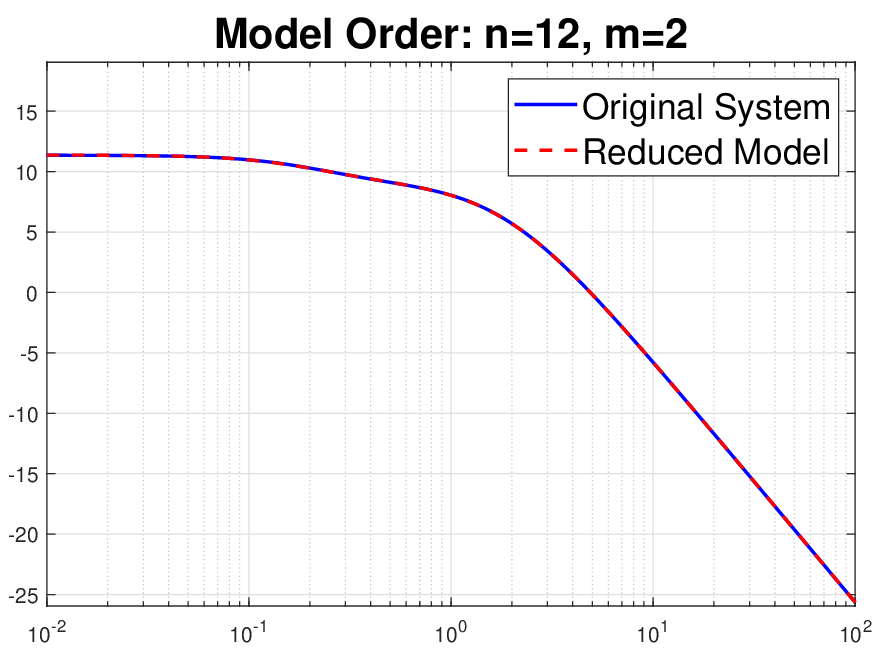}
    \end{subfigure}
    \hfill
    \begin{subfigure}{0.239\textwidth}
        \centering
        \includegraphics[width=\linewidth]{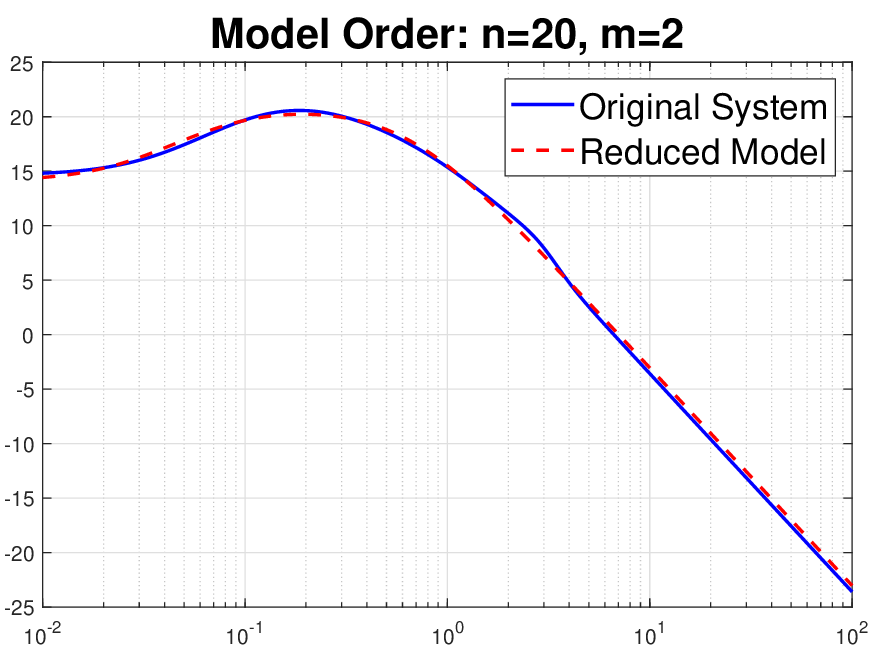}
    \end{subfigure}
    \begin{subfigure}{0.239\textwidth}
        \centering
        \includegraphics[width=\linewidth]{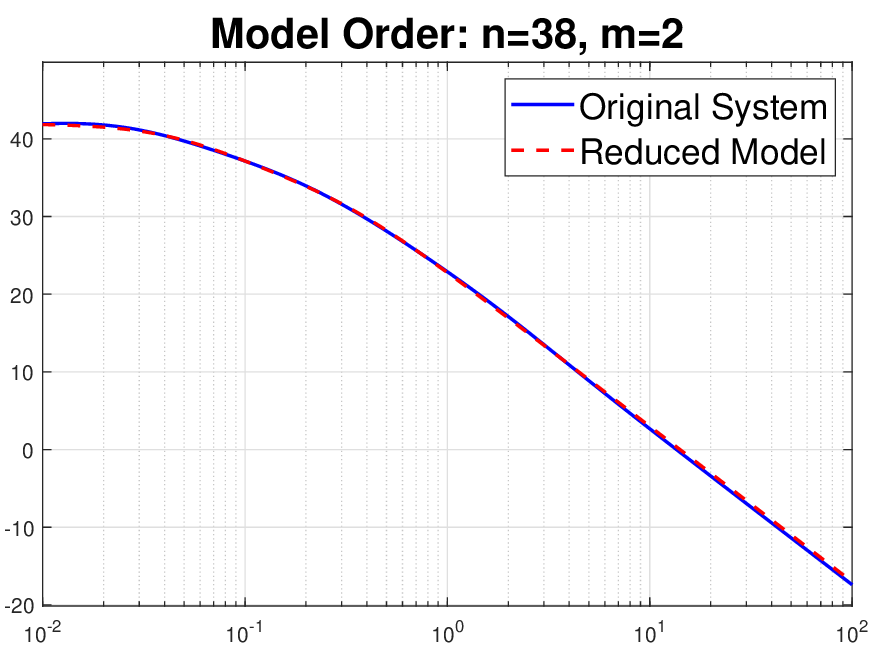}
    \end{subfigure}
    \hfill
    \begin{subfigure}{0.239\textwidth}
        \centering
        \includegraphics[width=\linewidth]{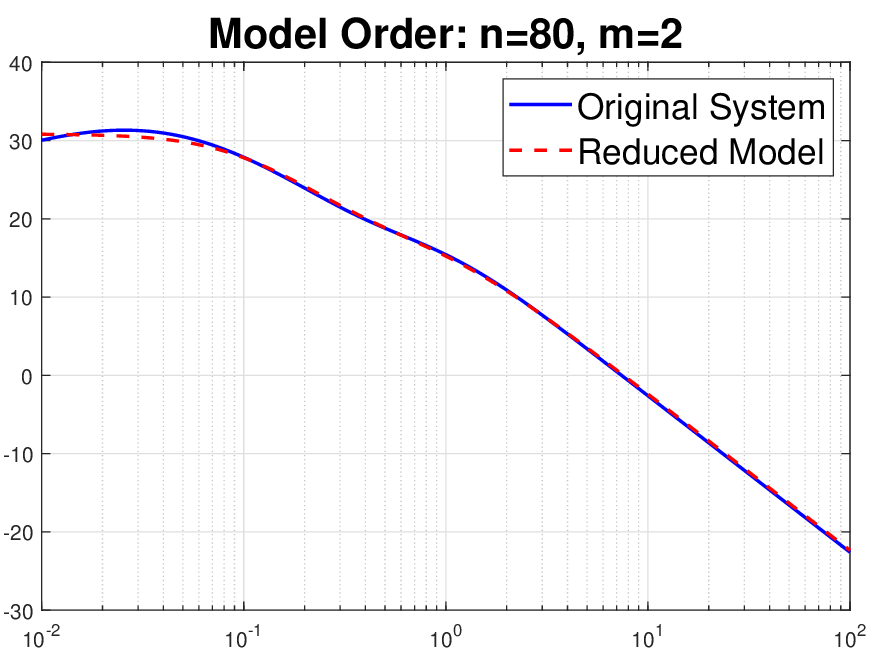}
    \end{subfigure}
\caption{Amplitude bode plots of the $G(s)$ ($n = 12, 20, 38, 80$) and $G_m(s)$ ($m = 2$) generated using the SDR method. The x-axis is the frequency in rad/s, while the y-axis denotes the magnitude in dB.}
\label{fig:combined}
\end{figure}

Further examples for $m=2$ are generated using the command \texttt{sys = rss(n,1,1)} for $n=12$ and $n=20$, with additional commands to avoid zero-pole cancellations. For $n=38$ and $n=80$, systems are generated using \texttt{randn}. The Bode diagrams in Fig.~\ref{fig:combined} demonstrate that the reduced-order models closely approximate the full models. 

\begin{figure}[htbp]
    \centering
    \includegraphics[width=\linewidth]{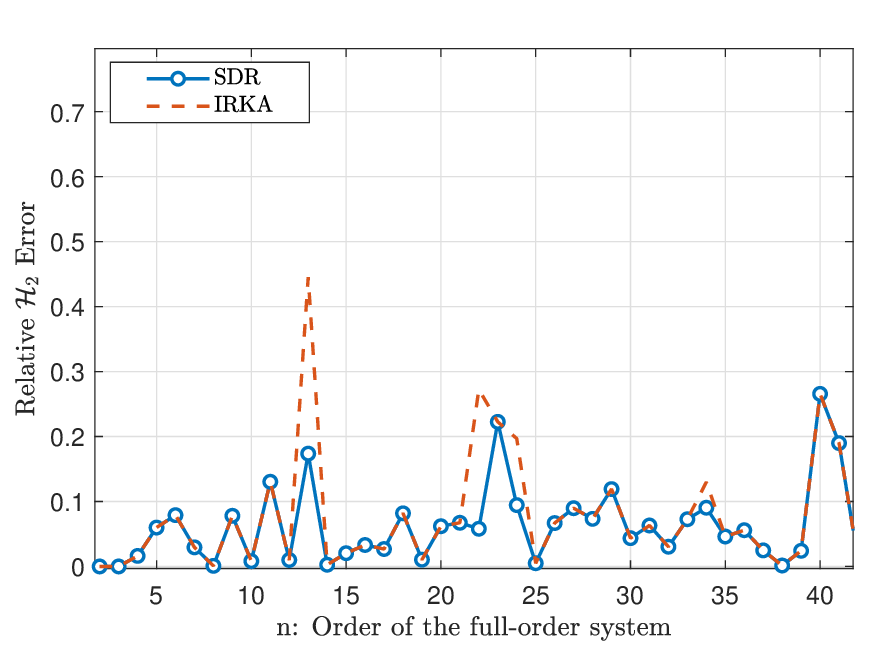}  
    \caption{Relative $\mathcal{H}_2$ norm of the error system vs. $n$.}
    \label{fig:relative error}
\end{figure}

Two different strategies were tested for generating a set of interpolation points to construct the reduced-order system: (1) the globally optimal shifts obtained via the proposed SDR method, and (2) randomly selected shifts from the interval $(-\infty, +\infty)$, which were then applied in IRKA. The results for $m = 2$, showing the relative $\mathcal{H}_2$ error for each $n$, are presented in Figure~\ref{fig:relative error}. The figure demonstrates that both strategies perform well. However, the optimal shifts from SDR consistently lead to lower errors than those obtained with randomly initialized shifts in IRKA for all values of $n$. We emphasize that these results are achieved through globally optimal shift selection, enabling the direct construction of reduced-order models. This makes the approach numerically effective for large-scale problems and allows it to overcome the global convergence challenges faced by other algorithms.
\vspace{-3mm}
\section{CONCLUSIONS}\label{sec:conclusion}
This paper addresses the $\mathcal{H}_2$ model reduction problem by developing an SDR  method that provides an upper bound on the optimal solution. The proposed SDR is formulated as LMIs and is applied for cases $m=1$ and $m=2$. This approach efficiently yields global optimal interpolation points that can be used in other iterative techniques to directly derive the reduced-order system. The exactness of the proposed upper bound is analyzed using numerical examples. We also give a proof of global optimality of our solution when $m = 1$. An extension of the SDR method to higher-order systems, will be reported in future work.

\appendix[Proof of Global Optimality for $m=1$]
\label{appendix}
\renewcommand{\theequation}{A.\arabic{equation}}
\setcounter{equation}{0}
Let $\mathbf{S}^n$ and $\mathbf{H}^n$ denote the sets of real symmetric and complex Hermitian $n \times n$ matrices, respectively. The generalized Kalman–Yakubovich–Popov (gKYP) lemma is concerned with inequalities on curves $\boldsymbol{\Lambda}(\Phi, \Psi)$ in the complex plane:
\begin{equation}\label{eq:lambda111}
\boldsymbol{\Lambda}(\Phi, \Psi):=\{\lambda \in \mathbb{C} \mid \sigma(\lambda, \Phi)=0, \sigma(\lambda, \Psi) \geq 0\},
\end{equation}
where $\Phi,\Psi \in \mathbf{H}^2$ and where for Hermitian $H$ we define
$$
\sigma(\lambda, H):=\left[\begin{array}{ll}
\lambda & 1
\end{array}\right]H\left[\begin{array}{l}
\bar{\lambda} \\
1
\end{array}\right]\!\cdot
$$
The next result defines $\Phi$ and $\Psi$ such that $\Lambda(\Phi, \Psi)$ represents the non-negative real axis.
\begin{lem}\label{lemma:restricted lambda}
Let $\Lambda(\Phi, \Psi)$ be defined as in (\ref{eq:lambda111}) with
\begin{equation}\label{eq:segment real axis}
\Phi=\left[\begin{array}{cc}
0 & i \\
-i & 0
\end{array}\right], \quad \Psi=\left[\begin{array}{ll}
0 & 1 \\
1 & 0
\end{array}\right]\cdot
\end{equation}
where $i$ denotes $\sqrt{-1}$. Then for all $\lambda \in \Lambda(\Phi, \Psi)$, $\lambda$ is restricted to be real and non-negative.
\end{lem}
\begin{proof}
Let $\lambda \in \mathbb{C}$. Then
$$
\sigma(\lambda, \Phi)=\left[\begin{array}{ll}
\lambda & 1
\end{array}\right]\left[\begin{array}{cc}
0 & i \\
-i & 0
\end{array}\right]\left[\begin{array}{l}
\bar{\lambda} \\
1
\end{array}\right]=i \lambda-i \bar{\lambda}=0
$$
capturing $\lambda \in \mathbb{R}$, and 
$$
\sigma(\lambda, \Psi)=\left[\begin{array}{ll}
\lambda & 1
\end{array}\right]\left[\begin{array}{ll}
0 & 1 \\
1 & 0
\end{array}\right]\left[\begin{array}{l}
\bar{\lambda}\\
1
\end{array}\right]=\lambda+\bar{\lambda} \geq 0
$$
capturing the real part of $\lambda \geq 0 $.
\end{proof}
\begin{lem}[generalized KYP Lemma\cite{kyp2,kyp5,kyp6}]\label{lem:gKYP}
Let $A \in \mathbb{C}^{n \times n}, B \in \mathbb{C}^{n \times m},\Theta \in \mathbf{H}^{n+m}$ and $\Phi, \Psi \in \mathbf{H}^2$ be given and define $\boldsymbol{\Lambda}(\Phi, \Psi)$ by (\ref{eq:lambda111}). Suppose that $(A,B)$ is controllable and that $A$ has no eigenvalues $\lambda$ such that $\sigma(\lambda, \Phi)=0$. Then, the following statements are equivalent:
\begin{itemize}
    \item[(i)] The inequality
    \vspace{-1mm}
    \begin{equation}
\left[\begin{array}{c}
(\lambda I-A)^{-1} B \\
I
\end{array}\right]^H \Theta\left[\begin{array}{c}
(\lambda I-A)^{-1} B \\
I
\end{array}\right] \preceq 0,
\vspace{-1mm}
\end{equation}
holds for all $\lambda \in \boldsymbol{\Lambda}(\Phi, \Psi)$.
    \item[(ii)] There exist Hermitian $P$ and $Q$  such that $Q \succ 0$ and
\begin{equation}
F^H(\Phi \otimes P+\Psi \otimes Q)F + \Theta \preceq0, \ \ F=\left[\begin{array}{cc}
A & B \\
I & 0
\end{array}\right] \text {. }
\end{equation}
\end{itemize}
\end{lem}
\vspace{1.5mm}
\begin{rem}\label{rem:Pure Imaginary}
The matrix $F$ can be replaced by any arbitrary matrix depending on the system. When $F,~\Theta,~\Phi$, and $\Psi$ are real, $P$ and $Q$ can be restricted to be real. This follows the fact that the real
part of a complex Hermitian positive–definite matrix is positive. As such, when $F,\Theta, \Psi$ are real, and $\Phi$ is purely imaginary, we may select $Q$ real and $P$ imaginary without loss of generality \cite{kyp3,kyp7}.  
\end{rem}
\begin{thm}\label{thm:KYP}
Let $T_4 \in \mathbb{R}^{n \times n}, C^T \in \mathbb{R}^{n \times m},\Theta \in \mathbf{S}^{n+m}$ and $\Phi, \Psi \in \mathbf{H}^2$ be given and define $\boldsymbol{\Lambda}(\Phi, \Psi)$ by (\ref{eq:lambda111}). Suppose $(T_4,C^T)$ is controllable and $T_4$ has no eigenvalues $\lambda$ such that $\sigma(\lambda, \Phi)=0$. Then, the following statements are equivalent.
    \begin{itemize}
    \item[(i)] The inequality
    \begin{equation*}
\left[\begin{array}{c}
\left(\lambda I-T_4\right)^{-T} C^T \\
I
\end{array}\right]^H \Theta\left[\begin{array}{c}
\left(\lambda I-T_4\right)^{-T} C^T \\
I
\end{array}\right] \preceq 0, 
\end{equation*}
 \text {holds for all } $\lambda \in \boldsymbol{\Lambda}(\Phi, \Psi)$.
    \item[(ii)] There exist Hermitian matrices $P$ and $Q \succ 0$ such that
    \begin{equation*}
\left[\begin{array}{cc}
T_4^T & C^T \\
I & 0
\end{array}\right]^H(\Phi \otimes P+\Psi \otimes Q)\left[\begin{array}{cc}
T_4^T & C^T \\
I & 0
\end{array}\right]+\Theta \preceq 0.
\end{equation*}
\end{itemize}
\end{thm}
\vspace{1.5mm}
\begin{proof}
The result follows from Lemma \ref{lem:gKYP} by selecting $F=\bigl[ \begin{smallmatrix}
T_4^T & C^T \\
I & 0
\end{smallmatrix} \bigr],$
which in turn implies that $\Theta$ is real.
\end{proof}
\vspace{1.5mm}
\begin{prop}\label{prop:exact proof}
Let all variables be as in Proposition~\ref{prop:m=1} and assume that $\gamma^2$ is minimized subject to the LMIs in \eqref{C_1}. Then
\begin{equation}\label{eq_g2}
\min_{\stackrel{G_1(s)~\rm{is~stable}}{\rm{dim}(G_1)=1}}\left\|G-G_1\right\|^2_{\mathcal{H}_2}=\gamma^2
\end{equation}
\vspace{-2mm}

\end{prop}
\begin{proof}
We prove the proposition by reformulating the LMI optimization in Proposition~\ref{prop:m=1} in the form of Theorem~\ref{thm:KYP}.

Starting from \eqref{eq:001}:
\begin{equation}
f(p_1)-\gamma^2=C \hat{\Delta} T_3+T_3^T \widehat{\Delta}^T C^T-C \widehat{\Delta} T_5 \widehat{\Delta}^T C^T-\gamma^2\leq 0,
\end{equation}

where $T_2$ is replaced by $-C$. Equivalently, we have
\begin{equation}\label{eq:equivalent}
\left[\begin{array}{ll}
C \hat{\Delta} & I
\end{array}\right]\left[\begin{array}{cc}
-T_5 & T_3 \\
T_3^T & -\gamma^2
\end{array}\right]\left[\begin{array}{c}
\hat{\Delta}^T C^T \\
I
\end{array}\right] \leq 0,
\end{equation}
where $\widehat{\Delta}=\Delta\left(I-T_4 \Delta\right)^{-1}$ and $\Delta=p_1 I$. Define $\lambda=\frac{1}{p_1}$. It follows that $\lambda\ge0$ since $p_1\ge0$. Therefore, it follows from Lemma \ref{lemma:restricted lambda} that $\lambda\in\Lambda(\Phi,\Psi)$ where $\Phi$ and $\Psi$ are defined in (\ref{eq:segment real axis}). An equivalent statement to (\ref{eq:equivalent}) is given by Theorem \ref{thm:KYP} (ii). Define $\Theta=\left[\begin{smallmatrix}
-T_5 & T_3 \\
T_3^T & -\gamma^2
\end{smallmatrix}\right]$, we have
\begin{equation}\label{eq:withoutS}
\left[\begin{array}{cc}
T_4^T&C^T \\
I&0
\end{array}\right]^H(\Phi \otimes P+\Psi \otimes Q)\left[\begin{array}{cc}
T_4^T&C^T \\
I&0
\end{array}\right]+\Theta \preceq 0 .
\end{equation}
Set
\begin{equation}
(\Phi \otimes P+\Psi \otimes Q)=\left[\begin{array}{cc}
0 & Q+i P \\
Q-i P & 0
\end{array}\right].
\end{equation}
Convert (\ref{eq:002}) to a negative form and extract the second term.
\begin{equation}
v \underbrace{\left[\begin{array}{cc}
-T_5+\left(S T_4^T+T_4 S^T\right) & T_3+S C^{\mathrm{T}} \\
T_3^T+C S^T & -\gamma^2
\end{array}\right]}_{\hat{L}}v^T
 \leq 0
\end{equation}
where $v = \left[\begin{array}{cc}
C \hat{\Delta} & 1
\end{array}\right]$. $\hat{L}$ is further expanded as
\begin{equation}\label{eq:withS}
\hat{L}\!=\!\!\left[\!\!\begin{array}{cc}
T_4^T & C^T \\
I & 0
\end{array}\!\!\right]^H\!\left[\!\!\begin{array}{cc}
0 & S^T \\
S & 0
\end{array}\!\!\right]\!\left[\!\!\begin{array}{cc}
T_4^T & C^T \\
I & 0
\end{array}\!\!\right]\!\!+\!\Theta \preceq 0.
\end{equation}
Note that there is an equality between (\ref{eq:withS}) and (\ref{eq:withoutS})  when $S\!=\!Q\!-\!iP$. Thereby, $S^T$ is replaced by $S^H$. It follows from Remark~\ref{rem:Pure Imaginary} that we can select $Q$ real and $P$ imaginary without loss of generality, resulting in $S\!\in \mathbb{R}^{n\times n}$ as we defined.
\end{proof}
\begin{rem}
When $m=1$, the optimal $\mathcal{H}_2$ model reduction problem can be solved by the convex program and there is no gap between the reformulated convex problem and the original model reduction problem.   
\end{rem}

\bibliographystyle{ieeetr}
\bibliography{main.bib} 
\end{document}